\documentclass[fleqn,11pt,twoside]{article}

\usepackage{amsthm,amsthm,amssymb, color, xcolor,epsfig, graphics, subfigure}

\usepackage{amsmath, graphicx, latexsym, lscape }

\makeatletter
\newcommand{\copyrightnote}[2]{{\renewcommand{\thefootnote}{}
 \footnotetext{\small\it
\begin{flushleft}
 \copyright \ #1   #2  
\end{flushleft}}}}

\newcommand{\Name}[1]{\begin{flushleft}
                       \LARGE \bf #1
                       \end{flushleft}\vspace{-3mm}}

\newcommand{\Author}[1]{\begin{flushleft}
                       \it #1 \end{flushleft}}

\newcommand{\Address}[1]{\begin{flushleft}
                       \it #1 \end{flushleft}}

\newcommand{\Date}[1]{\begin{flushleft}
                      \small  \it #1 \end{flushleft}}

\newcommand{\evenhead}{Author \ name}
\newcommand{\oddhead}{Article \ name}

\renewcommand{\@evenhead}{
\hspace*{-3pt}\raisebox{-15pt}[\headheight][0pt]{\vbox{\hbox to \textwidth
{\thepage \hfil \evenhead}\vskip4pt \hrule}}}
\renewcommand{\@oddhead}{
\hspace*{-3pt}\raisebox{-15pt}[\headheight][0pt]{\vbox{\hbox to \textwidth
{\oddhead \hfil \thepage}\vskip4pt\hrule}}}
\renewcommand{\@evenfoot}{}
\renewcommand{\@oddfoot}{}

\long\def\@makecaption#1#2{%
  \vskip\abovecaptionskip
  \sbox\@tempboxa{\small \textbf{#1.}\ \ #2}%
  \ifdim \wd\@tempboxa >\hsize
    {\small \textbf{#1.}\ \ #2}\par
  \else
    \global \@minipagefalse
    \hb@xt@\hsize{\hfil\box\@tempboxa\hfil}%
  \fi
  \vskip\belowcaptionskip}

\newcommand{\JNMPnumberwithin}[3][\arabic]{%
  \@ifundefined{c@#2}{\@nocounterr{#2}}{%
    \@ifundefined{c@#3}{\@nocnterr{#3}}{%
      \@addtoreset{#2}{#3}%
      \@xp\xdef\csname the#2\endcsname{%
        \@xp\@nx\csname the#3\endcsname .\@nx#1{#2}}}}%
}

\newcommand{\resetfootnoterule} {
  \renewcommand\footnoterule{%
  \kern-3\p@
  \hrule\@width.4\columnwidth
  \kern2.6\p@}
}

\renewcommand{\footnoterule}{}

\makeatother

\theoremstyle{definition}

\begin{document}

\renewcommand{\evenhead}{ {\LARGE\textcolor{blue!10!black!40!green}{{\sf \ \ \ ]ocnmp[}}}\strut\hfill 
C Rogers and A C Briozzo}
\renewcommand{\oddhead}{ {\LARGE\textcolor{blue!10!black!40!green}{{\sf ]ocnmp[}}}\ \ \ \ \   
Moving boundary problems with Ermakov symmetry reduction
}

%%%% Matter for the first page 
\thispagestyle{empty}
\newcommand{\FistPageHead}[3]{
\begin{flushleft}
\raisebox{8mm}[0pt][0pt]
{\footnotesize \sf
\parbox{150mm}{{Open Communications in Nonlinear Mathematical Physics}\ \  \ {\LARGE\textcolor{blue!10!black!40!green}{]ocnmp[}}
\ \ Vol.5 (2025) pp
#2\hfill {\sc #3}}}\vspace{-13mm}
\end{flushleft}}

\FistPageHead{1}{\pageref{firstpage}--\pageref{lastpage}}{ \ \ Article}

\strut\hfill

\strut\hfill

\copyrightnote{The author(s). Distributed under a Creative Commons Attribution 4.0 International License}

\Name{Moving boundary problems with Ermakov symmetry reduction: nonlinear superposition principle and reciprocal transformation applications}

\Author{
Colin Rogers $^{1}$, Adriana C. Briozzo $^{2}$}

\Address{{$^1$} School of Mathematics and Statistics, The University of New South Wales\\Sydney NSW 2052, Australia\\{$^2$} Depto. Matem\'atica, FCE, Univ. Austral, Paraguay 1950-CONICET\\S2000FZF Rosario, Argentina.\\Email: c.rogers@unsw.edu.au,  abriozzo@austral.edu.ar}

\Date{Received September 9, 2025; Accepted September 28, 2025}

\setcounter{equation}{0}

\begin{abstract}
\noindent 
Moving boundary problems of Stefan-type for a novel third order nonlinear evolution equation with temporal modulation are here shown to be amenable to exact Airy-type solution via a classical Ermakov equation with its admitted nonlinear superposition principle. Application of the latter together with a class of involutory transformations sets the original moving boundary problem in a wide class with temporal modulation. As an appendix, reciprocally associated exactly solvable moving boundary problems are derived.
\end{abstract}

\label{firstpage}

%%%% The Article text starts here

\section{Introduction}
A nonlinear equation as introduced by Ermakov \cite{ermakov} has subsequently proved to constitute the canonical base member of multi-component Ermakov-type systems with extensive physical applications in both nonlinear physics and continuum mechanics \cite{rogers1}. Importantly, the original Ermakov equation admits a now classical nonlinear superposition principle. The latter may, notably, be applied in the analysis of moving shoreline hydrodynamics with an underlying boundary \cite{rogers2}. Indeed, the classical single component Ermakov equation has diverse applications such as, \textit{inter alia}, in the nonlinear elastodynamics of boundary-loaded hyperelastic tubes \cite{rogers3, shahinpoor}, oceanographic eddy pulsrodon evolution \cite{rogers4}, magnetogasdynamics \cite{rogers5} and the analysis of rotating gas cloud phenomena \cite{rogers6}.

Here, moving boundary problems of Stefan-type for a novel $3^{rd}$ order nonlinear evolution equation with temporal modulation are shown to be amenable to exact solution via symmetry reduction in a classical Ermakov equation with its  nonlinear superposition principle.

\section{A Ermakov Symmetry reduction}A novel third order nonlinear evolution equation with temporal modulation is introduced here according to
\begin{equation}\label{ec}
    u_t+u_{xxx}+\lambda (t+a)^{\mu}u^{-4}u_x=0,\quad \lambda, \mu \in \mathbb{R}
    \end{equation}
with symmetry reduction to be derived via the ansatz
    \begin{equation}\label{trans}
        u=(t+a)^{m}\Psi\left(\frac{x}{(t+a)^{n}}\right).
    \end{equation}
Thus, insertion of the latter into \eqref{ec} yields
\begin{equation}
    m\Psi-n\xi\Psi'+(t+a)^{-3n+1}\Psi'''+\lambda(t+a)^{\mu-4m-n+1}\Psi^{-4}\Psi'=0
\end{equation}whence, $m=-1/3$, $n=1/3$ together with $\mu=-2$ so that
\begin{equation}\label{ecc}
    \Psi'''-(1/3)(\xi\Psi)'+\lambda\Psi^{-4}\Psi'=0, \qquad \xi=\frac{x}{(t+a)^{n}}.
\end{equation}
Integration of \eqref{ecc} leads to the symmetry reduction
\begin{equation}\label{zeta}
    \Psi''-(1/3)\xi\Psi-(\lambda/3)\Psi^{-3}=\zeta, \quad \zeta \in \mathbb{R}
\end{equation} whence, with $\zeta=0$ and the scaling $\xi=\epsilon z$ where $\epsilon^{3}=-\frac{3}{2}$ there results a classical Ermakov-type equation
\begin{equation}\label {reduction}
    \Psi_{zz}+(z/2)\Psi=k^{*}\Psi^{-3}, \quad k^{*} \in \mathbb{R}.
\end{equation} This admits the nonlinear superposition principle
\begin{equation}
    \Psi=\sqrt{c_1\Omega_1^{2}+2c_2\Omega_1\Omega_2+c_3\Omega_2^{2}}
\end{equation}
wherein $\Omega_1$, $\Omega_2$ constitute a pair of linearly independent solutions of the auxiliary linear equation
\begin{equation}
    \Omega_{zz}+(z/2)\Omega=0
\end{equation} with constants $k^{*}$ together with $c_i$, $i=1,2,3$ such that \begin{equation}
    c_1c_3-c_2^{2}=\frac{k^{*}}{W^{2}}
\end{equation}
where $W=\Omega_1\Omega_{2z}-\Omega_{1z}\Omega_2$ is the constant Wronskian of $\Omega_1,\Omega_2$.

The general nonlinear superposition principle can be derived via Lie group invariance as in \cite{rogers2}. Here, in view of the Ermakov symmetry reduction \eqref{reduction}
\begin{equation}
    \Psi=\sqrt{c_1A_i^{2}(-2^{1/3}\xi/\epsilon)+2c_2A_i(-2^{1/3}\xi/\epsilon)B_i(-2^{1/3}\xi/\epsilon)+c_3B_i^{2}(-2^{1/3}\xi/\epsilon)}
\end{equation}
wherein $A_i$ and $B_i$ are the Airy functions of the $1^{st}$ and $2^{nd}$ kind respectively. The $c_i$, $i=1,2,3$ together with $k^{*}$ constitute available parameters.

\section{A class of Moving Boundary Problems}

Here, a class of Stefan-type moving boundary problems for \eqref{ec} is introduced, namely
\begin{equation}\label{ec1}
    u_t+u_{xxx}+\lambda (t+a)^{-2}u^{-4}u_x=0,\quad 0<x<S(t), \quad t>0 
    \end{equation}
    \begin{equation}
        u_{xx}(S(t),t)-(\lambda/3)(t+a)^{-2}u^{-3}(S(t),t)=L_mS^{i}(t)\dot{S}(t)
    \end{equation}
    \begin{equation}
        u(S(t),t)=P_mS^{j}(t)
    \end{equation}
together with
\begin{equation}\label{cero}
        u_{xx}(0,t)-(\lambda/3)(t+a)^{-2}u^{-3}(0,t)=H_0(t+a)^{k}
    \end{equation}
and the initial condition $s(0)=S_0$.

In the sequel, the moving boundary $x=S(t)=\gamma(t+a)^{1/3}$ is adopted whence $S_0=\gamma a^{1/3}$.
\subsection*{Boundary conditions}

(I) \begin{equation*}
        u_{xx}(S(t),t)-(\lambda/3)(t+a)^{-2}u^{-3}(S(t),t)=L_mS^{i}(t)\dot{S}(t), \quad t>0
    \end{equation*}
Insertion of the symmetry ansatz \eqref{trans} into the preceding yields
\begin{equation}
    \Psi''(\gamma)-(\lambda/3)\Psi^{-3}(\gamma)=L_m S^{i}(t)\dot{S}(t)(t+a)=\frac{1}{3}L_m\gamma^{i+1}(t+a)^{\frac{i+1}{3}}
\end{equation}whence $i=-1$ together with 

\begin{equation}
   \Psi''(\gamma)-(\lambda/3)\Psi^{-3}(\gamma)=\frac{1}{3}L_m.
\end{equation}
(II)
\begin{equation*}
        u(S(t),t)=P_mS^{j}(t), \quad t>0
    \end{equation*}
Accordingly, 
\begin{equation}
        \Psi(\gamma)(t+a)^{-1/3}=P_m\gamma^{j}(t+a)^{j/3}, \quad t>0
    \end{equation}
so that $j=-1$ and
\begin{equation}
    \Psi(\gamma)=P_m\gamma^{-1}.
\end{equation}
(III) 

\begin{equation*}
        u_{xx}(0,t)-(\lambda/3)(t+a)^{-2}u^{-3}(0,t)=H_0(t+a)^{k}, \quad t>0.
    \end{equation*}
This yields
\begin{equation}
       \Psi''(0)-(\lambda/3)\Psi^{-3}(0)=H_0(t+a)^{k+1}
    \end{equation} so that $k=-1$ and $H_0$ is determined by
\begin{equation}
       \Psi''(0)-(\lambda/3)\Psi^{-3}(0)=H_0
    \end{equation}

It is important to record that the nonlinear evolution equation \eqref{ec} which can admit classical Ermakov symmetry reduction may be embedded in a wide class of temporally-modulated equations which inherit this property via application of transformations of the type $T^{*}$ given by
\begin{equation}
    dt^{*}=\rho^{-2}(t)dt, \quad x^{*}=x, \quad    u^{*}=\rho^{-1}(t)u, \qquad T^{*}
\end{equation} 
which augmented by the relation  $\rho^{*}=\rho^{-1}$ admit the involutory property $T^{**}=I$. Application of $T^{*}$ to the moving boundary problems determined by \eqref{ec1} - \eqref{cero} may be made to embed it in a wide class of temporally-modulated such Stefan-type problems which admit integrable classical Ermakov  symmetry reduction. Involutory transformations such as $T^{*}$ have their genesis in the autonomisation of Ermakov-type systems \cite{athorne}. In modern soliton theory, spatial analogues of $T^{*}$ have been applied in \cite{rogers7} to link modulated coupled systems to their canonical unmodulated counterparts. The key characteristic solitonic properties of the latter such as invariance under certain  Bäcklund type transformations \cite{rogers8} and being amenable to the inverse scattering procedure \cite{ablowitz1} are thereby  inherited by the associated spatially modulated systems.

A detailed analysis of moving boundary problems of generalised Stefan-type was initiated in \cite{rogers9} for the solitonic Harry-Dym equation \cite{vassiliou}. This was motivated by a remarkable solitonic connection made in \cite{vasconcelos} between the latter and the classical Saffman-Taylor problem with surface tension \cite{saffman}. In subsequent developments, nonlinear moving boundary problems of Stefan-type have been shown to be amenable to analytic solution for a range of canonical soliton equations via Painlevé II symmetry reduction \cite{rogers10, rogers11, rogers12, rogers13, rogers14, rogers15}.
Hybrid Ermakov-Painlevé II systems were originally derived in \cite{rogers16} in the context of symmetry reduction of multi-dimensional coupled nonlinear Schrödinger systems of Manakov-type. The canonical base single component Ermakov-Painlevé II equation was obtained therein in the analysis of certain transverse wave motions in a generalised Mooney-Rivlin hyperelastic material. This novel extension of the classical Ermakov equation has subsequently proved to have diverse physical applications, notably, in Korteweg capillarity theory \cite{rogers17}, cold plasma physics \cite{rogers18} and the analysis of Dirichlet-type boundary-value problems in the Nernst-Planck electrolytic system \cite{amster}. A connection to the classical Painlevé XXXIV equation was established in \cite{rogers19}. Ermakov-Painlevé II symmetry reduction of a novel class of extended mKdV equations with temporal modulation and applications to associated moving boundary problem of Stefan-type is currently in progress.
\begin{appendix}

\section*{Appendix: A reciprocal associated Moving Boundary Problem}

Reciprocal type transformations as originally introduced by Bateman \cite{bateman1} concerned the derivation of a class of novel invariance properties in homentropic two-dimensional gasdynamics as associated with admitted conservation laws. It was subsequently established in \cite{bateman2} that the reciprocal relations constitute a particular class of Bäcklund transformations. Invariance under such multi-parameter reciprocal-type Bäcklund transformations has recently been extended to relativistic gasdynamics systems in \cite{rogers20, rogers21, rogers22}. In \cite{rogers23, rogers24}, 1+1 dimensional reciprocal transformations have been applied to solve classes of moving boundary problems of Stefan-type for nonlinear evolution equations which model heat conduction in a range of metals as delimited by Storm in \cite{storm}. A reciprocal transformation was subsequently applied in \cite{rogers25} to such moving boundary problems in order to determine the conditions for the onset of melting due to applied boundary flux. In soil mechanics, reciprocal transformations have application to the analysis of Stefan-type moving boundary problems descriptive of percolation of liquids through porous media \cite{rogers26,rogers27}, It is remarked that a reciprocal transformation of another type has application in the context of 3+1-dimensional discontinuity wave theory \cite{donato}.

In modern soliton theory, reciprocal transformations associated with admitted conservation laws were originally introduced in \cite{kingston}. Therein conjugation was made with the classical Bianchi permutability theorem of pseudo-spherical surface theory as derived via invariance of the 1+1- dimensional sine-Gordon equation under a Bäcklund transformation and whereby multi-soliton solutions may be algorithmically generated in an iterative manner \cite{rogers8}. Reciprocal transformations were subsequently applied in \cite{rogers28} to link the canonical AKNS and WKI inverse scattering schemes of \cite{ablowitz2} and \cite{wadati1} respectively. Certain classes of 1+1-dimensional solitonic hierarchies are connected via reciprocal transformations as in \cite{rogers29} and \cite{rogers30}. Reciprocal transformations in 2+1- dimensions were introduced in \cite{rogers31} and later applied in \cite{oevel} to establish novel connections between members of a solitonic triad consisting of the canonical  Kadomtsev-Petviashvili, modified  Kadomtsev-Petviashvili as 2+1-dimensional Dym hierachies.

Reciprocal transformations have been applied to moving boundary problems for a wide range of solitonic equations in \cite{rogers10, rogers11, rogers12, rogers13, rogers14, rogers15}. Thereby, the reciprocally associated moving boundary problems inherent Painlevé II-type integrability. The procedure was originally introduced in \cite{rogers9} to determine such reciprocal associates for the 1+1-dimensional soliton Harry Dym equation. In \cite{rogers13}, a class of Stefan-type moving boundary problems for the modified Korteweg de Vries (mKdV) equation was solved via application of a Painlevé II similarity reduction which involved the classical Airy equation. A reciprocal transformation was applied to derive a novel class of exactly solvable moving boundary problems for the base Casimir member of the compacton hierarchy as set down in \cite{olver}. In addition, application of a class of involutory transformations $T^{*}$ with origin in the autonomisation of the Ermakov-Ray-Reid system as in \cite{athorne} was made to isolate novel solvable moving boundary problems for mKdV equations with Ermakov-type temporal modulation. In \cite{rogers14}, reciprocal transformations were applied to a class of Stefan-type problems for the Korteweg-deVries equation to generate, in turn, Airy-type exact solution to associated moving boundary problems both for the canonical nonlinear evolution of magma theory and a novel reciprocal associate of the KdV equation which incorporates a source term. The Gardner equation \cite{miura}  which subsumes the mKdV and KdV solitonic equations has proved to have diverse physical applications, notably in plasma physics \cite{ruderman}, optical lattice theory \cite{wadati2} and the analysis of nonlinear oceanic wave propagation phenomena \cite{grimshaw}. In a recent development it has been derived in a nonlinear elastodynamic context \cite{coclite}. It is  remarked that certain Ermakov-type connections with transverse wave propagation in Mooney–Rivlin hyperelastic materials have been detailed in \cite{destrade,rogers32}. Moving boundary problems of Stefan-type have been shown to be amenable to exact solution via Painlevé II symmetry reduction in \cite{rogers15} both for the solitonic Gardner equation and a novel reciprocal associate with a source term. Ermakov spatial modulation of the Gardner equation via involutory transformations was detailed.

\section*{
    Reciprocal Transformation: Application to Moving Boundary Problem \eqref{ec1} - \eqref{cero}.}

In the present context with regard to the nonlinear temporally modulated evolution equation 
 \begin{equation}\tag{A1}\label{a1}
     u_t+u_{xxx}+\lambda(t+a)^{-2}u^{-4}u_x=0,
 \end{equation}
 under the reciprocal transformation $R^{*}$ given by
 \begin{equation}\tag{A2}\label{a2}
     dx^*=u \;dx+[-u_{xx}+\frac{\lambda}{3}u^{-3}(t+a)^{-2}]\; dt, \qquad t^*=t, \quad u^*=\frac{1}{u}
     \end{equation}
     there results
   \begin{equation}\tag{A3}
     dx=u^*\;dx^* + \left[\frac{\partial}{\partial x^*}\left(\frac{1}{u^*}\frac{\partial}{\partial x^*}\left(\frac{1}{u^*}\right)\right)-\frac{\lambda}{3} u^{*4}(t^*+a)^{-2}\right]\; dt^*,
 \end{equation}with compatibility condition
 \begin{equation}\tag{A4}
    \frac{\partial u^{*}}{\partial t^{*}}=\frac{\partial}{\partial x^*}\left[\frac{\partial}{\partial x^*}\left(\frac{1}{u^*}\frac{\partial}{\partial x^*}\left(\frac{1}{u^*}\right)\right)-\frac{\lambda}{3} u^{*4}(t^*+a)^{-2}\right].
 \end{equation}
 By virtue of the reciprocal connection $R^{*}$ of the latter to \eqref{a1}, importantly it inherits admissible Airy-type symmetry reduction.

 \subsection*{Reciprocal Moving Boundary Problem Application}

 Under the reciprocal transformation $R^{*}$, the class of moving boundary problems \eqref{ec1} -\eqref{cero} becomes 
\[
 \frac{\partial u^{*}}{\partial t^{*}}=\frac{\partial}{\partial x^*}\left[\frac{\partial}{\partial x^*}\left(\frac{1}{u^*}\frac{\partial}{\partial x^*}\left(\frac{1}{u^*}\right)\right)-\frac{\lambda}{3} u^{*4}(t^*+a)^{-2}\right], \quad x^*\big|_{x=0}<x^*<x^*\big|_{x=S(t)}, \quad t^{*}>0,
 \]
 \begin{equation}\tag{A5}
\frac{\partial}{\partial x^*}\left(\frac{1}{u^*}\frac{\partial}{\partial x^*}\left(\frac{1}{u^*}\right)\right)\frac{1}{u^*}-\frac{\lambda}{3} u^{*4}(t^*+a)^{-2}=L_mS^{i}\dot{S},\quad \text{on} \quad x^{*}|_{x=S(t)},\quad t^{*}>0,
 \end{equation}
 \[\frac{1}{u^{*}}=P_mS^{j}(t),\quad \text{on} \quad x^{*}\big|_{x=S(t)},\quad t^{*}>0,\]
\[\frac{\partial}{\partial x^*}\left(\frac{1}{u^*}\frac{\partial}{\partial x^*}\left(\frac{1}{u^*}\right)\right)\frac{1}{u^*}-\frac{\lambda}{3} u^{*4}(t^*+a)^{-2}=H_0(t^*+a)^{k},\quad \text{on} \quad x^{*}|_{x=0},\quad t^{*}>0\]
wherein $S(t)=\gamma(t+a)^{1/3}.$

In the preceding 
 \begin{equation}\tag{A6}
 dx^{*}\big|_{x=0}=[-u_{xx}+(\lambda/3)u^{-3}(t+a)^{-2}]dt\big|_{x=0}=(t+a)^{-1}[-\Psi''+(\lambda/3)\Psi^{-3}]dt\big|_{x=0}=0
\end{equation}
by virtue of \eqref{zeta} with $\zeta=0$. Hence, $x^{*}\big|_{x=0}$ is constant. It is seen than it is required that $H_0=0$.

In addition, 
\begin{equation}\tag{A7}
 dx^{*}\big|_{x=S(t)}=dx^{*}\big|_{x=\gamma(t+a)^{1/3}}=[(t+a)^{-1/3}\Psi(\gamma) \dot{S}(t)-L_mS^{i}(t)\dot{S}(t)]dt
 \end{equation}
 \[
 =(t+a)^{-}[\Psi(\gamma)-L_m \gamma^{-1}]\frac{\gamma}{3}dt
\]
 whence 
 \begin{equation}\tag{A8}
 x^{*}|_{x=S(t)}=\frac{\gamma}{3}[\Psi(\gamma)-L_m \gamma^{-1}]\ln(t^{*}+a)=S^{*}(t)
\end{equation}
upto an additive constant and the reciprocal initial condition on the moving boundary becomes 
 \begin{equation}\tag{A9}
 S^{*}(0)=\frac{\gamma}{3}[\Psi(\gamma)-L_m\gamma^{-1}]\ln(a)
\end{equation}

In conclusion, it is remarked that in \cite{rogers13}, a reciprocal class of moving boundary problems associated with the mKdV equation has likewise been derived with logarithmic reciprocal boundary $x^{*}=S^{*}(t^{*})$. The mKdV equation and the solitonic hierarchy of which it is the base member may be obtained in a purely geometric manner in connection with the particular motion of curves in a plane \cite{rogers8}. Intrinsic geometric formulation of such a type has been applied in magnetohydrodynamics to uncover underlying geodesic characteristics in a normal congruence of surfaces \cite{rogers33}.

\end{appendix}

\subsection*{Acknowledgements}

The present work has been partially sponsored by the Project PIP No 11220220100532CO CONICET-UA and the Project O06-24CI1901 Universidad Austral, Rosario, Argentina

\label{lastpage}
\end{document}